\documentclass[12pt]{article}
\usepackage{amsmath,amssymb,amsfonts,amsthm,graphicx}
\usepackage[usenames,dvipsnames]{color}
\newcommand{\E}{{\epsilon}}

\newcommand{\G}{{\mathcal G}}
\newcommand{\be}{\begin{equation}}
\newcommand{\ee}{\end{equation}}

\newcommand{\old}[1]{}

\newcommand{\yellow}[1]{}

\newcommand{\T}{{\tau}}

\newcommand{\M}{{\cal M}}
\newcommand{\MM}{\mathbb{M}}
\newcommand{\Q}{{q_{_G}}}

\newcommand{\hh}{\hspace{.03cm}}

\newtheorem{theorem}{Theorem}

\begin{document}
\date{}

\title {Phases in Large Combinatorial Systems} \author{ Charles
  Radin\thanks{Department of Mathematics, University of Texas, Austin,
    TX 78712. E-mail: {\tt radin@math.utexas.edu}} }

\maketitle

\begin{abstract}
\noindent This is a status report on a companion subject to extremal
combinatorics, obtained by replacing extremality properties with
emergent structure, `phases'. We discuss phases, and phase transitions,
in large graphs and large permutations, motivating and using the
asymptotic formalisms of graphons for graphs and permutons for
permutations. Phase structure is shown to emerge using entropy and
large deviation techniques.
\end{abstract}

\vfill \eject
\section{Introduction}
\label{SEC:Intro}

Consider the structure of these two families of problems, both
fundamental to extremal combinatorics \cite{AK}:
\begin{enumerate}
\item Determine those simple graphs with a given density of
  subgraphs $A$ (say triangles) and with
  the highest possible density of subgraphs $B$ (say edges) 
\item Determine those permutations with a given density of
  patterns $A$ (say \hbox{1\hh2\hh3}) 
and with the highest possible density of patterns $B$
  (say \hbox{3\hh2\hh1}) 
\end{enumerate}

\begin{figure}[htbp]
\center{\includegraphics[width=3.5in]{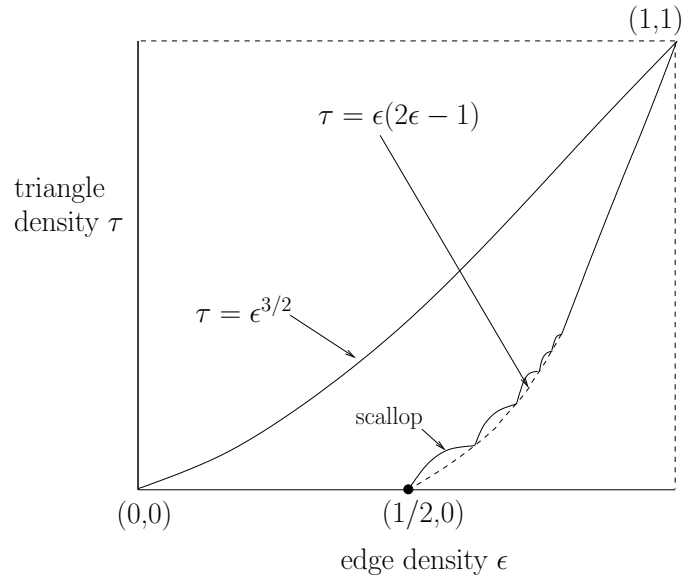}}
\caption{\label{scallop}Phase space for subgraph constraints edges and triangles}
\end{figure}

\begin{figure}[htbp]
\center{\includegraphics[width=3.5in]{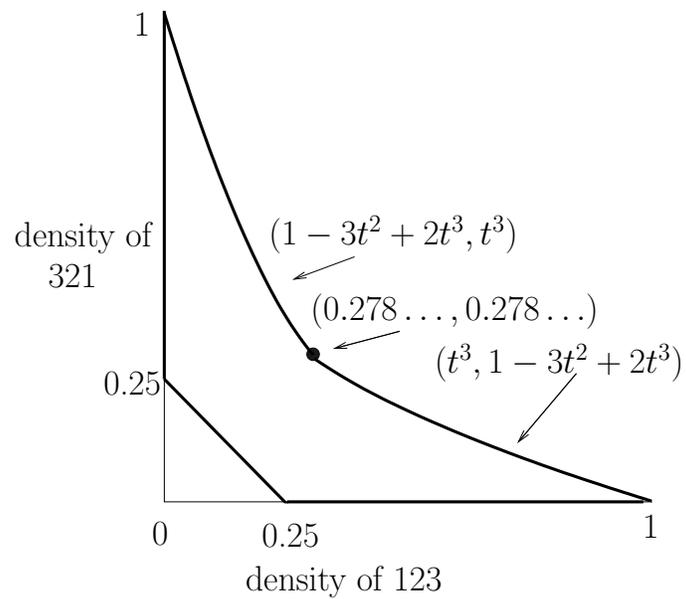}}
\caption{\label{123-321}Phase space for density constraints on the
patterns 1\hh2\hh3 and 3\hh2\hh1}
\end{figure}

\begin{figure}[htbp]
\center{\includegraphics[width=3.5in]{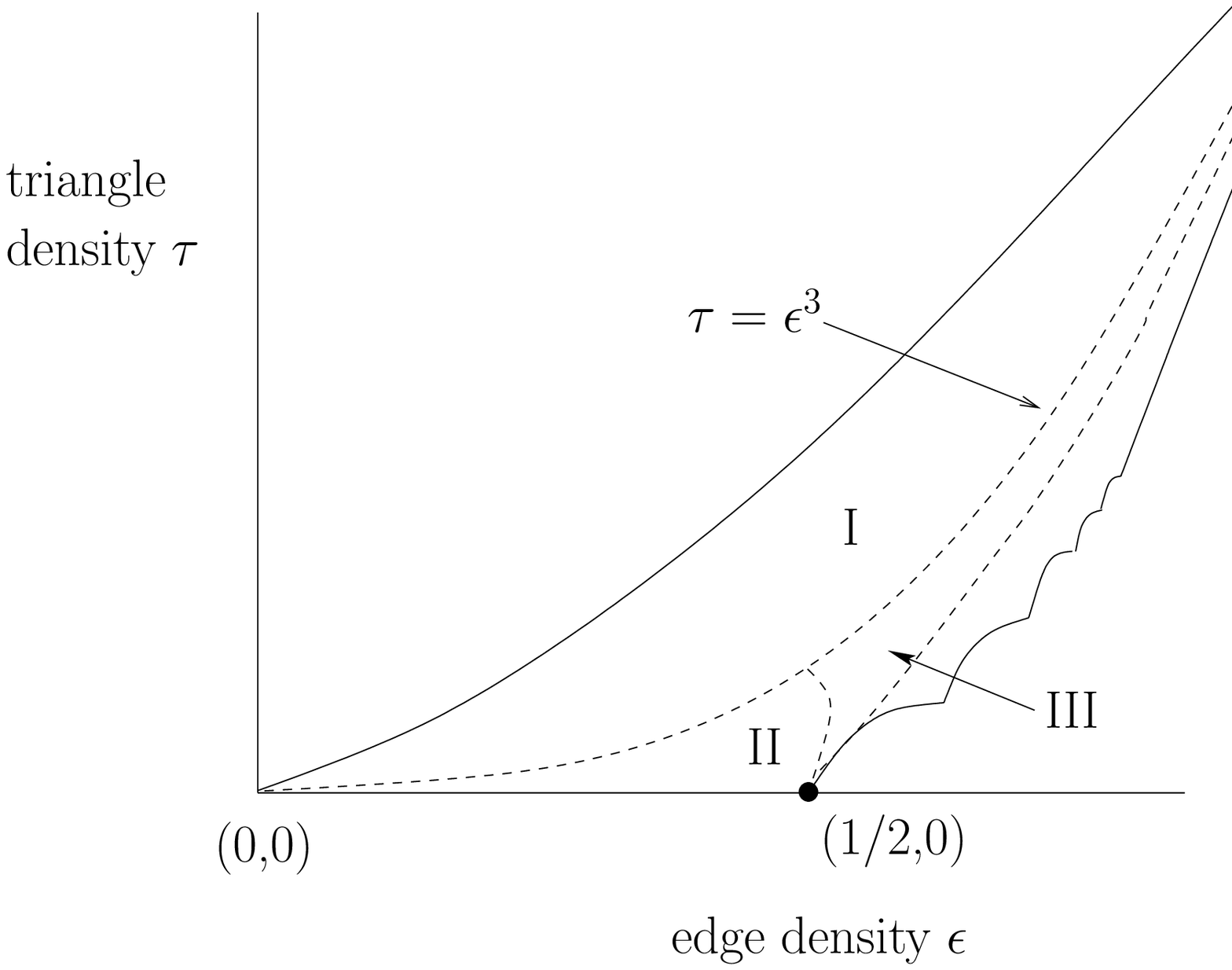}}
\caption{\label{scallop-phases}Three phases for edge/triangle constraints}
\end{figure}

\noindent For subgraphs, density refers to the fraction of them in a graph
compared to the number in the complete graph, while for a permutation
in $S_n$ a pattern is an element in $S_k,\ k\le n$, and the density of
the pattern in the permutation is the fraction of them compared to
${n\choose k}$.  We are only interested in asymptotic behavior, as
$n\to \infty$, where for instance $n$ is the number of nodes, for
graphs, or the number of objects $\{1,2,\cdots,n\}$ being permuted,
for permutations. For each such extremum problem we define the phase
space $\Gamma\subseteq [0,1]^2$ as the closure of the set of
simultaneously achievable pairs $(d_A, d_B)$ of densities of $A$ and
$B$, asymptotically in $n$.  See Figure
\ref{scallop} for an example of $\Gamma$ for graphs, and Figure
\ref{123-321} for an example for permutations.  (These and some other
figures are exaggerated to emphasize features.)
In this notation the above extremum problems 1 and 2 consist of determining
those graphs or permutations with density pairs on the boundary of
$\Gamma$. For some history of extremal graphs see \cite{Bo}; for
pattern avoidance in permutations see \cite{Ki}. In contrast to these
extremum problems, in this report we will be concerned with the
\emph{interior} of $\Gamma$, rather than its boundary, and not with
determining those graphs or permutations with such density
constraints, a hopeless and uninteresting problem, but with
determining what a \emph{typical one is like} for each achievable set
of constraints. (This is made precise in Section \ref{SEC: background}.)  The phenomenon on which
we focus is that well-defined phases emerge in $\Gamma$ as $n\to
\infty$, regions in $\Gamma$ in which typical behavior varies
smoothly with constraint values, separated however by sharp regional boundaries where
smoothness is lost (phase transitions), as indicated by dashed
curves in Figure \ref{scallop-phases}.

Analysis of such phases is greatly simplified by an asymptotic
formalism, associated with the term graphons for graphs (see \cite{L} for
a comprehensive introduction)
and permutons for permutations (developed in \cite{HKMS1,HKMS2}), and by large
deviation theorems in those asymptotic settings, as discussed in
\hbox{Section \ref{SEC: tools}.} 
The study of phases is more advanced for
graphs than permutations so this status report will mostly be about
graphs. There are very interesting extremum results in other parts of
combinatorics too, for instance partially ordered sets, but we do not
know of work on emergent phases in those fields (see however
\cite{PST} and references therein).
\vskip.4truein

\section{Introduction to phases in large contrained graphs}
\subsection{Background}
\label{SEC: background}
Before we wade in, here is some background.
The formalism of phases which we discuss in Section \ref{SEC: tools} mirrors that used in
statistical mechanics models, in which one analyzes configurations of $n$
particles in Euclidean space, with a specified potential energy
function whose gradient gives the interaction forces defining the
model. (For background in statistical
mechanics see \cite{Ru3}; there is a short outline in \cite{Ra4}.)
There is a long history studying configurations of
particles \emph{minimizing} the energy density for given mass
density (energy ground states) \cite{Ra1}, or \emph{maximizing} the mass density for
given energy density (densest packing) \cite{Fej, AW}. Solid and fluid phases emerge when
the achievable (energy density, mass density) pairs move away from 
the phase space boundary studied in these
extremum problems, and this is the phenomenon we are mirroring in
combinatorial systems. The physics models are much more complicated than the
combinatorial ones due to the geometric dependence
of the potential energy function, and indeed it is still an important open
problem \cite{Br, Uh, Si} to
prove the existence of a fluid/solid phase transition in any reasonably
satisfactory model (see however \cite{BLRW}). As we will see, progress has been quicker in
the simpler combinatorial settings. 

Our goal is to analyze `typical' large combinatorial systems with
variable constraints, for instance a typical large graph with
edge/triangle densities $(\E, \T)$ in the phase space of Figure
\ref{scallop}. Our densities are real numbers, limits of densities
which are attainable in large finite systems, so we begin by softening
the constraints, considering graphs with $n>>1$ nodes and with
edge/triangle densities $(\E^\prime,\T^\prime)$ satisfying
$\E-\delta<\E^\prime<\E+\delta$ and $\T-\delta<\T^\prime<\T+\delta$
for some small $\delta$ (which will eventually disappear.) It is easy to
show that the number of such constrained graphs is of the form $\exp
(sn^2)$, for some $s=s(\E,\T,\delta)>0$ and by a typical graph we mean
one chosen from the uniform  probability distribution $\mu_{n,\delta}^{\E,\T}$ on the constrained set.
\vskip.2truein 
\noindent
{\bf Goal.} \emph{We wish to analyze such families of constrained,
  uniform distributions of large combinatorial systems, in
  particular their dependence on the constraints.}
\vskip.2truein
\subsection{Main results for graphs}
\label{SUBSEC: Main}
Specializing to graphs, we present here the main qualitative results
on the asymptotics of constrained systems. Proofs are
scattered throughout \cite{RS1, RS2, RRS, KRRS1, KRRS2}, under various
hypotheses. The phenomena have also been verified to hold under weaker
restrictions by careful simulations, and we postpone describing these
issues of the range of validity until Section \ref{SEC: graph
  constraints}.
\begin{enumerate}
\item For fixed constraint values, and asymptotically in $n$, the
  nodes fall into a finite number of equivalence classes. More specifically,
in a typical large constrained graph the set of nodes can be partitioned into a finite
({\it usually small}) number of subsets, with well-defined fractions
of edges connecting nodes within, and between, such subsets. Thus
there are a small number of parameters describing our target
distributions, and they are functions of the constraints. For instance
using edge and triangle constraints, $\E$ and $\T$, for most values
of $(\E,\T)$ and asymptotically in $n$, the uniform
distribution has four parameters $a,b,c,d$: the node set is
partitioned into two subsets of relative size $c$, with edges between
a fraction $a$ of node pairs in one set, $b$ for the other set, and
fraction $d$ of node pairs split between the sets, and these are
functions of $(\E,\T)$.
\item The phase space $\Gamma$ of achievable constraint values can be
  partitioned into regions called phases, within which the
  parameters of the distribution are unique and smooth functions of
  the constraints, separated by lower dimensional boundaries on which
  parameter values may or may not be unique, but at which they lose
  their smoothness (phase transitions). Figure
  \ref{scallop-phases} exhibits three of the phases in the system with
  edge/triangle constraints. Phase $II$ exhibits an unusual level of symmetry, between
  {\em classes} of nodes, rather than merely within classes.
\end{enumerate}
Together these show how the phase structure is exhibited by typical
graphs, through the development of levels of equivalence among nodes.

\section{The basic tools: entropy and graphons}
\label{SEC: tools}
\subsection{Graphons}
\label{SEC: Graphons}
The mathematics of asymptotically large
graphs uses graphons, which we now review \cite{L}. The set of graphs, on $n$ 
nodes labelled $\{1,\ldots,n\}$, 
will be denoted $G_n$. (Graphs are assumed simple, i.e.\ undirected
and without
multiple edges or loops.  Given $G$ in $G_n$ we use
its adjacency matrix to represent $G$ by the function
$\Q$ on the unit square $[0,1]^2$ with constant value 0 or 1 in each
of the subsquares of area $1/n^2$ centered at the points $([j-1/2]/n,
[k-1/2]/n)$. More generally, a graphon $q\in \G$ is an
arbitrary symmetric
measurable function $[0,1]^2$ with values in $[0,1]$.
Informally, $q(x,y)$ is the probability of an edge between nodes $x$
and $y$, and so two graphons
are called equivalent if they agree up to a 
`node rearrangement' (see \cite{L} for details). 
Taking representatives, 
we define the cut metric on the quotient space $\tilde \G$ of `reduced
graphons' by
\be 
{d}(\tilde f,\tilde g)\equiv \inf_{f\in \tilde f, g\in \tilde g}\sup_{S,T\subseteq [0,1]}\Big| \int_{S\times
  T}[f(x,y)-g(x,y)]\, dxdy\Big|. 
\ee
\noindent $\tilde \G$ is compact in this topology \cite{L}. (We will define an
equivalent metric on $\tilde \G$ in \eqref{EQ: bar metric}.)

We now consider `blowing up' a graph $G$ by replacing each
node with a cluster of $K$ nodes, for some fixed $K=2,3,\ldots$, with
edges inherited as follows: there is an edge between a node in
cluster $V$ (which replaced the node $v$ of $G$) and a node in cluster
$W$ (which replaced node $w$ of $G$) if and only if there is an edge
between $v$ and $w$ in $G$. (The resultant graph is multipartite.)
Note that the blowups of a graph are all represented by the same
reduced graphon, and $\Q$ can therefore be considered a graph on
arbitrarily many -- even infinitely many -- nodes, which we exploit next.

The features of a graph $G$ on which we have focused are the
densities with which various subgraphs $H$ sit in $G$. Assume for
instance that $H$ is a quadrilateral. We could represent the density of
$H$ in $G$ in terms of the adjacency matrix $A_G$ by
\be
\frac{1}{{n\choose 4}} \sum_{w,x,y,z} A_G(w,x)A_G(x,y)A_G(y,z)A_G(z,w),
\ee

\noindent where the sum is over distinct nodes $\{w,x,y,z\}$ of $G$. 
For large $n$ this can approximated, within $O(1/n)$, as:

\be
\int_{[0,1]^4} \Q(w,x)\Q(x,y)\Q(y,z)\Q(z,w)\, dw\, dx\, dy\, dz.
\ee

\noindent It is therefore useful to define the
density $t_H(q)$ of this $H$ in a graphon $q$ by

\be \label{average}
\int_{[0,1]^4} q(w,x)q(x,y)q(y,z)q(z,w)\, dw\, dx\, dy\, dz,
\ee

\noindent and use such densities for this $H$, and analogs for other subgraphs, in
analyzing constrained distributions.  We note that $t_H(q)$ is a
continuous function of $q$ if we use the cut metric on reduced graphons. This
set of functions is also separating: any two reduced graphons with the
same values for all densities $t_H$ are the same  \cite{L}.

Next we give a different view of graphons. Let $\M$ be the set
of multipodal graphons, i.e. those for which there is a partition of
$[0,1]$ into finitely many subsets $F_j$ and the graphon is constant
on each product $F_j\times F_k$, and let $\M^\prime$ be the subset of
those functions of the form $\Q$ for some graph $G$.
Consider the metric $\bar d$ on reduced graphons defined by
\be \label{EQ: bar metric} \bar d(\tilde f,\tilde g)= \sum_{j\ge 1}
|t_{H_j}(\tilde f)-t_{H_j}(\tilde g)|/2^j, \ee 
\noindent where
$\{H_j\}$ is any ordering of the countable set of finite simple
connected graphs. This metric is equivalent on $\tilde \G$ to the cut
metric defined earlier \cite{L}. We can thus realize $\tilde \G$ in an
obvious way as a space of sequences, with  coordinates in $[0,1]$, and
metric $\bar d$, and note that the image $\tilde \G^\prime$ of $\M^\prime$
is dense in $\tilde \G$ \cite{L}.

\subsection{Entropy and the variational principle for graphs}
\label{SEC: entropy for graphs}
Getting back to our goal of analyzing constrained uniform
distributions, $\mu_{n,\delta}^{\alpha}$,
a related step is to determine the cardinality of the
set of graphs on $n$ vertices subject to constraints. Our constraints
are expressed in terms of a vector $\alpha$ of values of a set $C$ of densities, 
and a softening agent
$\delta$. Denoting the cardinality by $Z_n(\alpha,\delta)$, it was proven
in \cite{RS1, RS2} that $\lim_{\delta\to 0}\lim_{n\to \infty}
(1/n^2)\ln[Z_n(\alpha,\delta)]$ exists; it is called the constrained
entropy $s_{\alpha}$. From thermodynamics it is known that much can
be learned simply from knowing this function of the values $\alpha$ of the
constraints. In statistical mechanics one focuses differently,
using probabilistic notions to analyze the asymptotic
constrained uniform distributions, again as a function of the
constraint values, and this is what we discuss for these combinatorial
settings in this report.

There are $2^{n\choose 2}\approx e^{kn^2}$ graphs in $G_n$.  There
are at most $n!\approx \exp^{nln(n) -n}$ graphs equivalent to any
particular element of $G_n$, which for large $n$ is negligible
compared to $e^{kn^2}$ and this fact will be relevant below.

Suppose we are interested in analyzing those $G\in G_n$ with edge density
approximately $e_0\in (0,0.5]$ and the largest possible triangle
density, which is $(e_0)^{3/2}$ \cite{L}. To attain this one must use $m$
nodes to form a clique (all possible edges), where $m$ is determined
by $2^{m\choose 2}=2^{n\choose 3}e_0^{3/2}$, and leave the remaining
nodes as spectators (no connections). There are many ways to do
this, but they are all represented by the same reduced graphon in
$\tilde \G^\prime$.

Alternatively, suppose we wanted the $G$'s to have edge density $e_0\in
(0,0.5]$ but with minimal possible triangle density, which is 0. To
achieve this one can select two subsets $A,B$ of the nodes, and
choose $ne_0$ edges but only between nodes in different sets. There
are many inequivalent bipartite graphs of this type (except for
$e_0=0.5$), so a more productive goal might be to get a useful
handle on the distribution of solutions $G$.

Finally, there is an enormous number of ways to attain edge density $e_0\in (0.5]$ and
triangle density fixed between the maximum and minimum just
discussed. For these intermediate cases
we change the problem as suggested in the previous paragraph; we
no longer try to identify the appropriate graphs, but it turns out
we can often identify what a typical such graph is like, i.e.\
determine the (uniform) distribution on such constrained graphs.
For instance for triangle density $\tau\in [0,e_0^3]$ the constrained
distribution is obtained by partitioning the node set into two equal
sets $A$ and $B$, and choosing edges between node pairs independently
by the following two rules: for any two nodes in the same set the
probability of an edge is $[e_0-(e_0^3-\tau)^{1/3}]$,
and for any two nodes in different sets the probability is $[e_0
  +(e_0^3-t)^{1/3}]$. For $\tau=e_0^3$ this reduces to the situation
in which each pair of nodes is connected with probability $e_0$, while
for
$\tau=0$ it reduces to the extreme case discussed above.

This is the point where the mathematics changes
flavor. This is not due merely to our focus on asymptotics; as noted in Section
\ref{SEC:Intro}, the
extremal combinatorics associated with the boundary of the phase space
already involves asymptotics, and for instance led to the beautiful
flag algebra formalism of Razborov \cite{R3}. The difference here is
that at points in the interior of $\Gamma$, where we want to
understand not individual graphs but the constrained uniform
distribution on graphs, our problem is naturally reformulated
within the {\em calculus of variations}, since a key tool is the 
constrained entropy
$s_\alpha$ which can be represented as follows.
\begin{theorem}\label{thm:g-variation}
(The variational principle for constrained 
graphs \cite{RS1,RS2})
{For any vector $H$ of subgraphs $H_j$ and vector $\alpha$ of
  numbers $\alpha_j$,
\be \label{var}
s_\alpha=\max_{t_H(q)=\alpha} S(q),
\ee
\noindent  where
$S$ is the Shannon entropy of graphons:}
\be \label{Shannon}
S(q)=-\int_{[0,1]^2} \frac{1}{2}\{q(x,y)\ln[q(x,y)]+[1-q(x,y)]\ln[1-q(x,y)]\}\,dxdy.
\ee
\end{theorem}
Our goal is to understand families of constrained uniform
distributions, $\mu_{n,\delta}^{\alpha}$, 
on graphs on $n$ nodes and constraints $\alpha,\delta$,
with $n$ large and $\delta$ small. It can be tricky to
`compare' distributions for different $n$; we overcome this using
graph blowup as in Section \ref{SEC: Graphons} to work with graphons,
giving us a uniform framework independent of $n$. We then get our
approximation to $\mu_{n,\delta}^{\alpha}$ using the optimal
graphons $\tilde g^\alpha$ of Theorem \ref{thm:g-variation}. Any
reduced graphon $\tilde g$ is a sequence of densities $t_{H_j}$, and we
use the fact 
that the average values of these densities, with respect to the sequence
$\mu_{n,\delta}^{\alpha}$, converge to those of $\tilde g^\alpha$
(assuming $\tilde g^\alpha$ is a unique optimizer) \cite{RS1}. These limiting
averages are computable from $\tilde g^\alpha$ as in (\ref{average}),
and are our main use of the distributions $\mu_{n,\delta}^{\alpha}$.
(They would also allow us to directly estimate the
mass functions of the $\mu_{n,\delta}^{\alpha}$ if desired.) 
We can therefore interpret the averages as the densities of a `typical' large, constrained
graph. 

Variational principles such as Theorem \ref{thm:g-variation} are well
known in statistical mechanics \cite{Ru1,Ru2,Ru3}, and for their
simpler (discrete) models can be obtained from general large
deviations results \cite{El}. This was also the case for Theorem
\ref{thm:g-variation}, for which the proof used the large deviation
theory from \cite{CV}.  
One aspect of the applicability of such variational principles in this
context is not well understood. In all known graph examples the
optimizers in the variational principle are unique in $\tilde \G$
except occasionally on a lower dimensional set of constraints where
there is a phase transition. But we do not yet have a theoretical
understanding of this fundamental issue, sometimes called the Gibbs
phase rule in physics; see \cite{Ru2,I} for weak versions in physics. Given such
uniqueness however, the rest of the path is clear: the optimizer for
(\ref{var}) gives us the limiting constrained uniform distribution
which we want to analyze!  (Without uniqueness it is harder to obtain
useful information.)

\section{Examples of constrained graph systems}
\label{SEC: graph constraints}
We now give some details, including references, concerning the
qualitative results claimed in Section \ref{SUBSEC: Main}. As noted
there, our understanding of the range of validity of what is proven in
our theorems is significantly enhanced by careful simulations, in a
range of examples (models). To preserve continuity of argument we will
clarify some statements with references to the Notes at the end of the
text.

\subsection{Edge/triangle model}
\label{SEC: e/t model}
\begin{figure}[htbp]
\center{\includegraphics[width=3in]{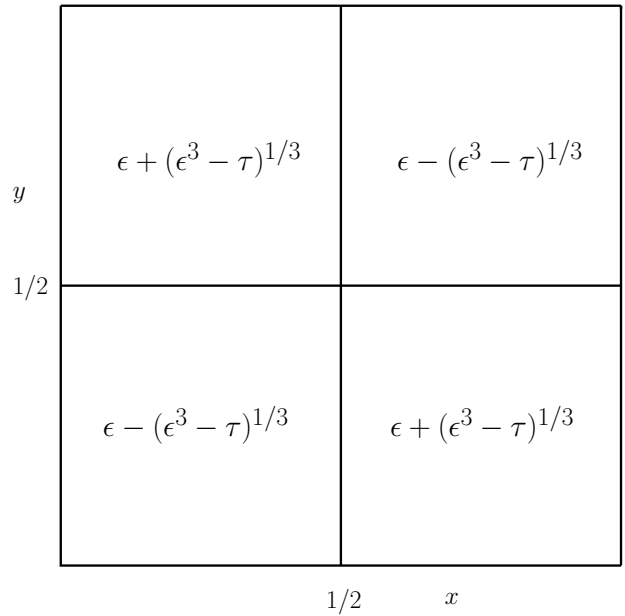}}
\caption{\label{graphon}The (piecewise constant) edge/triangle graphon in phase $II$}
\end{figure}

Consider first the edge/triangle model, using edge/triangle densities
$(\E,\T)$ as constraints, studied in
\cite{RS1,RS2,RRS,KRRS1,KRRS2}. The entropy optimizer is unique$^{(1)}$ in
$\tilde \G$ for
every pair $(\E,\T)$. Also, the optimizer
is multipodal, i.e.\ there is a partition of
$[0,1]$ into finitely many subsets $F_j$ and the graphon is constant
on each product $F_j\times F_k$. This gives a distribution with finitely many parameters, which vary
smoothly with $(\E,\T)$ except on certain curves, separating
phases$^{(2)}$. Specifically, for over $95\%$ of the area of
$\Gamma$, corresponding to constraints in the phases $I,\ II$ and
$III$ in Figure \ref{scallop-phases}, it is bipodal -- the partition of the nodes is into
just two subsets -- so there are only four parameters in the asymptotic
distribution, each a function of the constraints. 
In particular, in
phase $II$ there is a very interesting symmetry: the two {\em sets} of nodes are in all ways
interchangeable. This means, in terms of the parameters defined earlier,
that $a=b$ and $c=1$, so there are only two independent variables,
$a$ and $d$, in this phase. One can then solve for $a$ and $b$ in
terms of the constraints $\E$ and $\T$, getting the graphon in Figure
\ref{graphon} (using convenient node labelling). 

\begin{figure}[htbp]
\center{\includegraphics[width=3.5in]{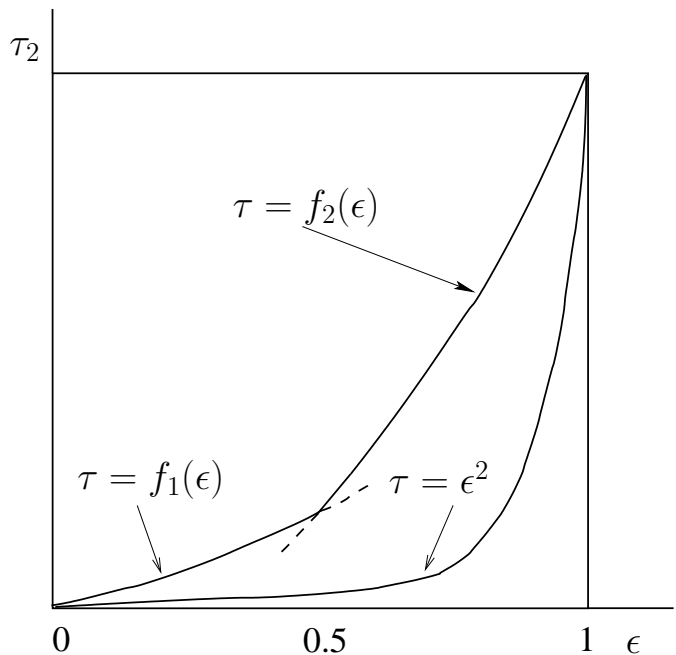}}
\caption{\label{2-star}The phase space for the 2-star model}
\end{figure}

\begin{figure}[htbp]
\center{\includegraphics[width=3.5in]{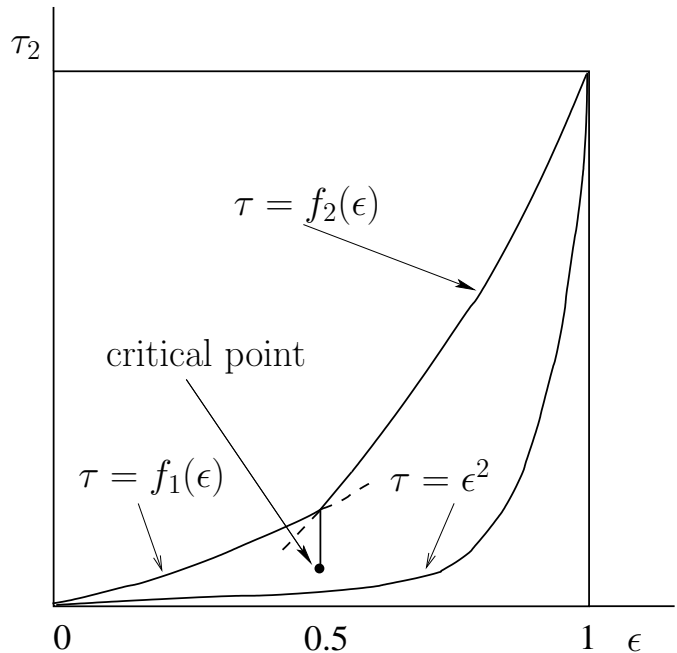}}
\caption{\label{critical}The critical point in the 2-star model}
\end{figure}

It is natural to study the boundaries of phases, to see how the system
changes there. Consider the approach to the phase boundary $\T=\E^3$, from each
side. Fixing $\E \le 0.5$ and increasing $\T$ from 0, we see from
Figure \ref{graphon} that the large typical graph starts as perfectly symmetric
bipartite, and steadily loses the distinction between the two sets$^{(3)}$. At
$\T=\E^3$ the nodes are all equivalent; these are the
Erd\H{o}s-R\'enyi graphs, with iid edges, and since the constrained
entropy peaks  at the curve $\T=\E^3$ when varying either $\E$ or $\T$
separately, these graphs are important in part because of their overwhelming
number. They also play the same intuitive role for us as the ideal gas
in statistical mechanics.

The approach to the Erd\H{o}s-R\'enyi curve from phase $I$ is quite
different. Except at $\E=0.5$ (the maximum of the unconstrained
entropy), as one increases $\T$ from $\E^3$ the optimal graphon
immediately becomes bipodal but now by the emergence of an
infinitesmal set of nodes: $c=O(\T-\E^3)$. The asymptotic form of
$a,b$ and $d$ are proven in \cite{KRRS2}. The behavior of phase $I$ near its upper
boundary, $\T=\E^{3/2}$ has not been analyzed. But on the upper
boundary itself the optimal graphon consists of spectator nodes, fully
unconnected, and a cluster of fully interconnected nodes (a clique),
of size needed to produce edge density $\E$. It is not yet understood
how these graphs connect to the ones emerging from the lower boundary
of phase $I$, though extensive simulation data has been gathered \cite{RRS}. And
finally, perturbation about the boundary between phases $III$ and $II$ is
particularly interesting because of the breaking of symmetry.
New features appear in other models, which we now discuss. 

\subsection{$k$-star graph models}
A $k$-star is a graph with $k$ edges with a node in common, and the
\hbox{$k$-star model} is the one in which the constraints are the
density $\E$ of edges and the density $\T_k$ of $k$-stars. The phase
space $\Gamma$ of the 2-star model is typical of them, and shown in
Figure \ref{2-star} \cite{KRRS1}$^{(4)}$. For any $k$-star model the
lower boundary of $\Gamma$ is the curve, $\T_k=\E^k$, represented by
Erd\H{o}s-R\'enyi (constant) graphons, and the upper boundary is the
upper envelop of two intersecting smooth curves$^{(4)}$. There is a
unique bipodal graphon representing all but one point on the upper
boundary, the intersection point just noted, where the two one-sided
graphon limits differ$^{(4)}$. We know for $k$-star models that in the
interior of $\Gamma$ there is a unique entropy-optimizing bipodal
graphon everywhere except on a curve emanating from the intersection
point on the boundary; this curve ends at a `critical point' in the
interior of $\Gamma$, so there is only one phase in each of these
models$^{(5)}$. See Figure \ref{critical}. The behavior of $k$-star
models just above the Erd\H{o}s-R\'enyi curve is common to many
other models, including the edge/triangle model \cite{KRRS2}. However
this universality does not extend further from the Erd\H{o}s-R\'enyi
curve; in particular there is no transition in the edge/triangle model
analogous to that in the $k$-star models$^{(6)}$. (There is some
confusion on this point coming from exponential random graph
models$^{(7)}$.)

\subsection{Half-blip model}\label{SUBSEC: half-blip}
Finally we consider the half-blip model, where the constraints are
a pair of signed densities, the signed 2-star density $t_1$ and signed
square \hbox{density $t_2$,} defined for graphons by:
\be 
t_1(q)= \int q(x,y) [1-q(y,z)] dx\, dy\, dz ;
\ee

\be
t_2(q)=\int q(w,x)
[1-q(x,y)] q(y,z) [1-q(z,w)] dw \, dx \, dy \, dz.
\ee

\begin{figure}[htbp]
\center{\includegraphics[width=2.5in]{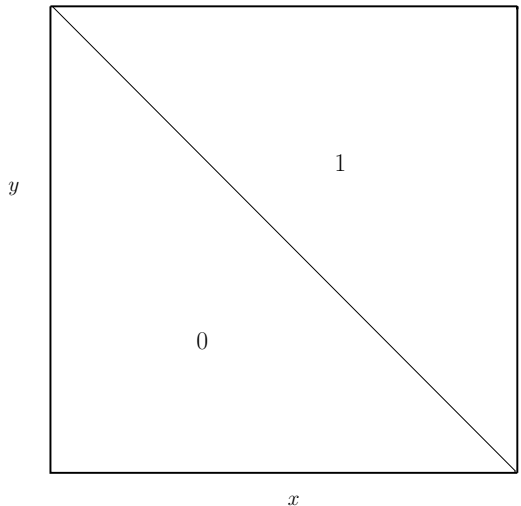}}
\caption{\label{half-blip}The half-blip graphon at $(t_1,t_2)=(1/6,0)$}
\end{figure}

\noindent The phase space for this model is not fully known, but there is a
lower edge corresponding to $t_2=0$, namely $\{(t_1,0)\,|\, 0\le
t_1\le 1/6\}$, with the following feature: as $t_1$ increases from 0
the unique representing graphon is $m$-podal but with $m\to\infty$ as
$t_1\to 1/6$. This feature is not so special, as it is shared with the
edge/triangle model along its lower (scalloped) boundary as $\E\to
1$. However in the half-blip model the unique graphon associated with
$(t_1,t_2)=(1/6,0)$ is not multipodal -- it is the graphon, with
values 0 and 1, shown in Figure \ref{half-blip} -- and this unusual circumstance is
quite intriguing.  

There is much to learn about this model, both on
the boundary and the interior of $\Gamma$. The model is an
analog of statistical mechanics models of quasicrystals \cite{Ra2}, which
raises the question whether, in the half-blip model, there is a phase
in the interior of $\Gamma$ near $(1/6,0)$ without multipodal entropy
optimizers. Indeed the analogous question has not yet been solved in
statistical mechanics \cite{AR}, but might well be easier to solve in this
combinatorial setting.

We have discussed several models of large graphs, but all had two
constraints. We do not know of new phenomena which become available
with more constraints, though we note one model studied with three
constraints, the taco model \cite{KRRS1}. In the other direction, a model with
only one constraining density reduces to the Erd\H{o}s-R\'enyi family,
which is very important but does not exhibit phases, or transitions, in
the senses we have used the terms. There has been much work in
adjusting the Erd\H{o}s-R\'enyi models, for instance to exhibit
percolation. This is of course important, but seems to
be exploring phenomena essentially different, say, from the emergence of
the symmetric phase $II$ in large graphs with edge/triangle constraints, as
described in Section \ref{SEC: e/t model}.

\section{Phases in large constrained permutations}
\subsection{Permutons and entropy}
Next we review some results in which the above style of
analysis of phases is applied to large constrained permutations.  We
begin with a quick review of the asymptotic permuton formalism, introduced in \cite{HKMS1,HKMS2}.
Given a labelled set of $n$ objects, which we take to be $\{1,2,\ldots,n\}$,
a permutation
$\pi$ in $S_n$ is an invertible function from $\{1,2,\ldots,n\}$ to
itself, represented  
by its values $\displaystyle (\pi_1,\pi_2,\ldots,\pi_n).$
It is commonly displayed as a square $0-1$
matrix $A_{j,k}$ with value 1 when $\pi_j=k$.  This can be rescaled
and reinterpreted as a function on $[0,1]^2$ which has values 0 or 1
on each square $[(j-1)/n, j/n]\times [(k-1)/n, k/n]$,
$j,k=1,2,\ldots,n$, with value 1 when $\pi_j=k$. Since such functions
are nonnegative and have integral 1, they can be reinterpreted as
the densities of probability measures $\gamma_\pi$ on $[0,1]^2$, with the property of
uniform marginals:
\be
\gamma_\pi([a,b]\times [0,1])=b-a=\gamma_\pi([0,1]\times [a,b]),
\hbox{ for all }
0\le a\le b\le 1.
\ee

\noindent Permutons are then defined more generally to be probability
measures on $[0,1]^2$ with uniform marginals. We put the usual weak
topology of measure theory (or weak$^\ast$ topology from functional
analysis) on the compact space $\MM$ of permutons, and note that the set
$\cup_n\{\gamma_\pi\,|\,\pi\in S_n\}$ is dense in $\MM$ \cite{HKMS1,HKMS2,GGKK,KKRW}.

We will be considering asymptotic conditional distributions on $S_n$,
asymptotic as $n\to \infty$, and we choose the conditioning to mesh with
the (extremal) study of pattern avoidance as indicated in the beginning of this
report. This choice is a more serious decision in the study of permutations
than conditioning by subgraph densities for graphs, because it emphasizes
a {\em linear} ordering of the objects being permuted,
which is quite restrictive; with another choice one might for instance
use the permutations to model the mixing of objects in space. In any
case, as we did with graphs we define the density $\rho_\tau(\gamma)$
of a pattern $\tau\in S_k$ in a permuton $\gamma$ as an asymptotic
form of its natural meaning for permutations, namely by the
probability that when $k$ points are selected independently from
$\gamma$ and their $x$-coordinates are ordered, the permutation
induced by their $y$-coordinates is $\tau$.  For example, for $\gamma$
with probability density $g(x,y)\,dxdy$, the density of pattern $12\in
S_2$ in $\gamma$ is
\be
\rho_{12}(\gamma)=2\int_{x_1<x_2\in[0,1]}\int_{y_1<y_2\in[0,1]}
g(x_1,y_1)g(x_2,y_2)\,dx_1dy_1dx_2dy_2.
\ee

\noindent  It follows that if $\gamma_{\pi_j}$ converges to $\gamma$ then
the density of $\tau$ in $\pi_j$ converges to $\rho_\tau(\gamma)$,
and that two permutons are equal if they
have the same pattern densities for all patterns. See \cite{GGKK} for background on
permutons. 

We now use pattern densities to condition permutations.
Let $\gamma$ be a permuton with probability density $g$.
We define the {Shannon entropy} $H(\gamma)$ of $\gamma$ by:
\be
H(\gamma) = \int_{[0,1]^2} -g(x,y) \ln g(x,y)\,dxdy,
\ee
where
$0 \ln 0$ is taken as 0.  Then $H$ is finite whenever $g$ is
bounded (and sometimes when it is not).  In particular for any
$\pi \in S_k$, we have $H(\gamma_\pi)= -k(k \ln{k} /k^2) = -\ln k$
and therefore $H(\gamma_\pi) \to -\infty$ for any sequence of
increasingly large permutations even though $H(\lim \gamma_\pi)$ may
be finite.  Note that $H$ is 0 on the uniform permuton (where
$g(x,y) \equiv 1$) and negative (sometimes $-\infty$) on all other
permutons, since the function $-z \ln z$ is concave downward.  If
$\gamma$ has no probability density we define $H(\gamma) = -\infty$.

Fix some finite set $\{\tau_1,\dots,\tau_k\}$ of patterns
and let $\alpha = (\alpha_1,\dots,\alpha_k)$ be a vector of their target
densities. We then define two sets of permutons:
\be
\Lambda^{\alpha,\delta}=\{\gamma\in \MM\,|\, |\rho_{\tau_j}(\gamma)
- \alpha_j| < \delta \hbox{ for each }1 \le j \le k\},
\ee
\be
\Lambda^{\alpha}=\{\gamma\in \MM\,|\, \rho_{\pi_j}(\gamma) = \alpha_j 
\hbox{ for each }1 \le j \le k\}.
\ee
With that notation, and the understanding that
$\Lambda_n^{\alpha,\delta} = \Lambda^{\alpha,\delta} \cap
\gamma(S_n)$, where $\gamma(\pi)=\gamma_\pi$, we have:
\begin{theorem}\label{thm:p-variation}
(The variational principle for constrained 
permutations \cite{KKRW})
{\be
\lim_{\delta\to 0}\lim_{n\to\infty}\frac{1}{n} \ln
\frac{|\Lambda_n^{\alpha,\delta}|}{n!}
= \max_{\gamma \in \Lambda^{\alpha}} H(\gamma).\ee}
\end{theorem}
Constrained sets of permutations in $S_n$ have cardinality of
order $\displaystyle e^{n \ln n +(c-n)}$ where $c \in [-\infty,0]$ is
our target \cite{KKRW}.
The function of $\alpha$, $\displaystyle \max_{\gamma \in \Lambda^{\alpha}} H(\gamma)$, which is
guaranteed by the theorem to exist but may be $-\infty$, is the 
\emph{constrained entropy} and denoted by $s(\alpha)$. Theorem \ref{thm:p-variation} was proven
in \cite{KKRW} using the large deviations theorem from \cite{Tr}.

\subsection{Examples of constrained large permutations}
The general framework for the asymptotics of constrained permutations
is thus analogous to that of constrained graphs. In detail the mathematics
is quite different, but in terms of results about phase structure 
the main difference is in the depth of progress on
examples, even on the boundary of their phase spaces, i.e. the extremal
theory. The main examples in which relevant progress has been made concerning
phases in constrained permutations are:
the 1\hh2/1\hh2\hh3 model, star models, and the 1\hh2\hh3/3\hh2\hh1 model.

For the 1\hh2/1\hh2\hh3 model it has been shown \cite{KKRW} that the phase space
is the same scalloped triangle which is the phase space for
the edge/triangle of graphs \hbox{Figure \ref{scallop}.} For graphs the vertices of the
scalloped triangle were shown to give rise to interesting phases and
phase transitions, as indicated in \hbox{Figure \ref{scallop-phases}.} However there is no evidence
yet of an analogous phenomenon for the \hbox{1\hh2/1\hh2\hh3} model of permutations.

Star models of permutations use constraints of slightly different
character than considered so far. Instead of a single pattern $\tau\in
S_k$ say, one replaces two or more of the symbols in $\tau$ by a $*$,
which is a place holder which can be filled by any unused symbols.
For instance the constraint $*2*$ fixes the density of the union of
all consistent patterns in $S_3$, namely 1\hh2\hh3 and 3\hh2\hh1. For a class of
models with a finite number of such constraints it is proven that
constrained entropy is an analytic function of the constraint values,
and that it has a unique optimal permuton at each point in its phase
space \cite{KKRW}. 

Finally, for the 1\hh2\hh3/3\hh2\hh1 model there is proven \cite{KKRW}
to be a phase
transition on a curve emanating from the singularity at
$(0.278\ldots,0.278\ldots)$ on its phase boundary (see Figure \ref{123-321}), of
similar character to the one proven for the \hbox{edge/2-star} graph model
(see Figure \ref{2-star}).

We conclude our discussion of the permutation theory by noting two
topics in the asymptotics of constrained permutations which were
omitted due to a lack of results on phase structure. One is the
useful tool of insertion measures; see \cite{KKRW}. And finally, we
have ignored work done with a single constraint because, as was noted
in Section \ref{SUBSEC: half-blip} for graphs, it does not seem to bear on the deeper issues of
phases and phase transitions; however see \cite{SW} and references
therein.

\section{Open problems}

Two open problems for constrained graphs go to the heart of our
understanding of phase emergence. The first is to understand why, in
all models we have studied, there is a unique optimizer for the
constrained entropy except off a set of constraints of lower
dimension. As noted in Section \ref{SEC: entropy for graphs}, without
this feature our method would not produce useful results. There is
some analysis of this phenomenon in statistical mechanics, but the
only results are uniqueness off a set of constraint values of category
one, or of measure zero \cite{Ru2,I}, and while intuitively
suggestive this is not of practical use.

Another basic problem is the origin of the multipodal structure of
entropy optimizers. We have never found a nonmultipodal entropy
optimizer in the interior of a phase space, i.e.\ in a phase, but
there is no general proof that they cannot exist. (There are proofs in
some regions of some models \cite{KRRS1,KRRS2}, 
but they do not give any insight into the general situation.) An obvious candidate
for exploring this is the half-blip model of Section \ref{SUBSEC:
  half-blip}, a model which had already been used to study the
analogous phenomenon in extremal graph theory \cite{LS}.

\begin{figure}[htbp]
\center{\includegraphics[width=3.5in]{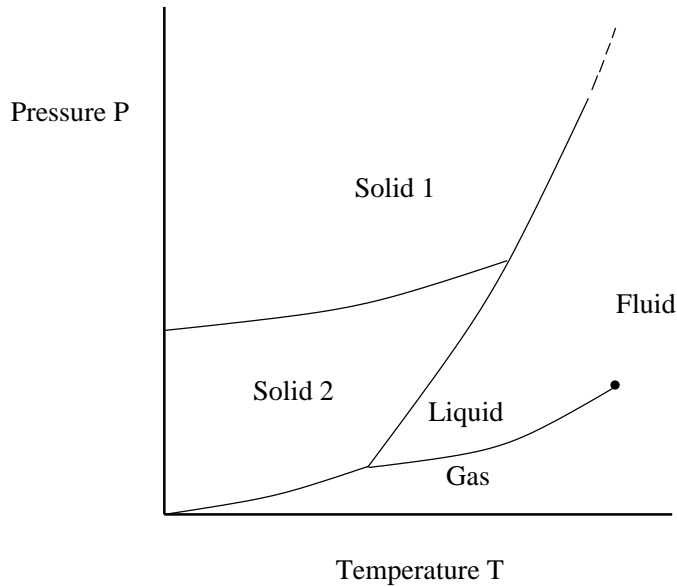}}
\caption{\label{PT diagram}The phase diagram of a simple material}
\end{figure}

Multipodality can be understood in terms of a symmetry between nodes,
as discussed in Section \ref{SUBSEC: Main}.  Another, and related, open
problem concerns phase $II$ in the edge/triangle
model. The symmetry of this phase, discussed in Section \ref{SEC: e/t
  model}, is between the equivalence classes that define multipodal
structure. What is the significance of this higher level of symmetry?
Multipodal graphons are piecewise constant as functions on $[0,1]^2$,
and it is plausible that such uniformity originates from
maximizing entropy, as was actually demonstrated in \cite{KRRS1} for $k$-star models.
 But the higher level of symmetry of phase $II$ in the edge/triangle
 model is
in a different category. One possible path to understanding this
links the problem  to an old open problem in statistical physics, to understand
why there is not a critical point in the transition between solid and
fluid phases of matter as there is between liquid and gas; see Figure
\ref{PT diagram}. 

Consider the traditional symmetry argument from physics \cite{An}:
\vskip.1truein \noindent
``It was Landau (Landau and Lifshitz, 1958) who, long ago, first pointed
out the vital importance of symmetry in phase transitions. This, the
First Theorem of solid-state physics, can be stated very simply: it is
impossible to change symmetry gradually. A given symmetry element is
either there or it is not; there is no way for it to grow
imperceptibly. This means, for instance, that there can be no critical
point for the melting curve as there is for the boiling point: it will
never be possible to go continuously through some high-pressure phase
from liquid to solid.''
\vskip.1truein 
This argument is not completely convincing \cite{Pi}, but the role of
symmetry in the study of solids is pervasive and has been highly
productive both in science and mathematics. And in support for linking 
the problems we note there is no critical
point in the transition between phases $II$ and $III$ in Figure
\ref{scallop-phases} where the graphons for phase $II$ but not $III$ have the higher
symmetry, while there is a critical point in Figure
\ref{critical}, where the graphons do not exhibit the higher level symmetry.

Concerning permutations we note that the permuton machinery seems to
apply well to pattern constraints. But this leaves out a major aspect
of permutations, their cycle structure, which makes use of permutation
multiplication and avoids the linear structure used in
patterns. Incorporating cycle structure should be a major source of development.

\vfill \eject
{\centerline{\bf Notes}}
\begin{enumerate}

\item This has been proven in a thin region above the
Erd\H{o}s-R\'enyi curve \cite{KRRS2} and on the line $\{(\E,\T)\,|\,\E=0.5,\
0\le \T\le 1/8\}$ \cite{RS1}, and seen in all the extensive simulations of the
model noted in \cite{RRS}.

\item Bipodality is proven in a thin region above the
Erd\H{o}s-R\'enyi curve  \cite{KRRS2} and on the line segment $\{(\E,\T)\,|\,\E=0.5,\
0\le \T\le 1/8\}$ \cite{RS1}. Bipodality and tripodality has been seen in all
other regions simulated in \cite{RRS}. There are expected to be regions of
$m$-podality with increasing $m$ near the corners of the scallops
along the boundary of the phase space \cite{RS1,RS2}. The rest of the claims
about entropy optimizers in this paragraph are only known from the
simulations referenced.

\item The bipodal graphon in Figure \ref{graphon} is only proven along the
line segment $\{(\E,\T)\,|\,\E=0.5,\ 0\le \T\le 1/8\}$ \cite{RS1}, but has been
thoroughly investigated by simulation \cite{RRS} in the region described as
phase $II$.

\item The phase diagrams for $k$-star models are derived, but only for $k\le 30$, 
in \cite{KRRS1}.

\item The one phase is smooth except on the curve. Nonsmoothness on
  the curve has been
  proven for $k=2$, and seen in simulation in $k=3$, in \cite{KRRS1}.

\item This is only known from simulation, but would not be
expected because of the behavior of the edge/triangle model on the
upper boundary of its phase space.

\item Exponential random graph models are widely used, especially
in the social sciences, to model graphs on a fixed, small number of
nodes \cite{N}. These models are sometimes considered Legendre transforms of the
models being discussed in this report \cite{RS1}. However, as was pointed out in
\cite{CD}, as the number of nodes gets large the parameters in the model
become redundant, and this confuses any interpretation of phase
transitions in such models. One way to understand the difficulty
is that the constrained entropy in these models is not convex or
concave, and the Legendre transform is not invertible \cite{RS1}.
\end{enumerate}

\section*{Acknowledgments}

This work was partially supported by NSF grants DMS-1208191 and
DMS-1509088.

\end{document}